\documentclass[10pt, titlepage]{article}

\usepackage{amssymb}
\usepackage[polish, english]{babel}
\usepackage[pdftex]{color, graphicx}
\usepackage{float}
\usepackage{rotating}
\usepackage{authblk}
\usepackage{amsmath}
\usepackage{appendix}

\begin{document}

\title{The FEM approach to the 2D Poisson equation in 'meshes' optimized with the Metropolis algorithm}
\author{I.D. Kosi\'nska}
\selectlanguage{polish}
\affil{Wroclaw University of Technology \\
Institute of Biomedical Engineering and Instrumentation \\
Wybrze\.ze Wyspia\'nskiego 27, 50-370 Wroc"law, Poland}
\date{}
\maketitle

\selectlanguage{english}

\begin{abstract}
The presented article contains a 2D mesh generation routine optimized with the Metropolis algorithm. The procedure enables to produce meshes with a prescribed size $h$ of elements. These finite element meshes can serve as standard discrete patterns for the Finite Element Method (FEM). Appropriate meshes together with the FEM approach constitute an effective tool to deal with differential problems. Thus, having them both one can solve the 2D Poisson problem. It can be done for different domains being either of a regular (circle, square) or of a non--regular type. The proposed routine is even capable to deal with non--convex shapes.
\end{abstract}

\section{Introduction}
The variety of problems in physics or engineering is formulated by appropriate differential equations with some boundary conditions imposed on the desired unknown function or the set of functions. There exists a large literature which demonstrates numerical accuracy of the finite element method to deal with such issues \cite{zienkiewicz}. Clough seems to be the first who introduced \emph{the finite elements} to standard computational procedures \cite{clough}. A further historical development and present--day concepts of finite element analysis are widely described in references \cite{zienkiewicz,zienkiewicz2}. In Sec. 2 of the paper and in its Appendixes \ref{A}-\ref{D}, the mathematical concept of the Finite Element Method is presented.\\ 
In presented article the well-known Laplace and Poisson equations will be examined by means of the finite element method applied to \emph{an appropriate 'mesh'}. The class of physical situations in which we meet these equations is really broad. Let's recall such problems like heat conduction, seepage through porous media, irrotational flow of ideal fluids, distribution of electrical or magnetic potential, torsion of prismatic shafts, lubrication of pad bearings and others \cite{laplace}. Therefore, in physics and engineering arises a need of some computational methods that allow to solve accurately such a large variety of physical situations.\\
The considered method completes the above-mentioned task. Particularly, it refers to a standard discrete pattern allowing to find an approximate solution to continuum problem. At the beginning, the continuum domain is discretized by dividing it into a finite number of elements which properties must be determined from an analysis of the physical problem (e. g. as a result of experiments). These studies on particular problem allow to construct so--called \emph{the stiffness matrix} for each element that, for instance, in elasticity comprising material properties like stress -- strain relationships \cite{clough,clough2}. Then the corresponding \emph{nodal loads}\footnote{Nodes are mainly situated on the boundaries of elements, however, can also be present in their interior.} associated with elements must be found.\\
The construction of accurate elements constitutes the subject of a mesh generation recipe proposed by the author within the presented article. In many realistic situations, mesh generation is a time--consuming and error--prone process because of various levels of geometrical complexity. Over the years, there were developed both semi--automatic and fully automatic mesh generators obtained, respectively, by using the mapping methods or, on the contrary, algorithms based on the Delaunay triangulation method \cite{delaunay_lit}, the advancing front method \cite{afm} and tree methods \cite{tm}. It is worth mentioning that the first attempt to create fully automatic mesh generator capable to produce valid finite element meshes over arbitrary domains has been made by Zienkiewicz and Phillips \cite{phillips}.\\ 
The advancing front method (AFM) starts from an \emph{initial} node distribution formed on a basis of the domain boundary, and proceeds through a sequential creation of elements within the domain until its whole region is completely covered by them. The presented mesh algorithm takes advantage from the AFM method as it is demonstrated in Sec. 3. After a node generation along the domain boundary (Sec. 3.1), in next steps interior of the domain is discretized by adding \emph{internal} nodes that are generated at the same time together with corresponding elements which is similar to Peraire \emph{et al.} methodology \cite{Peraire}, however, positions of these new nodes are chosen differently according to the manner described in Sec. 3.2. Further steps improve the quality of the mesh by applying the Delaunay criterion to triangular elements (Appendix \ref{E}) and by a node shifting based on the Metropolis rule (Sec. 4).

\section{The finite element method}

\subsection{The mathematical concept of FEM}

\begin{figure}[t!]
\includegraphics[scale=0.3]{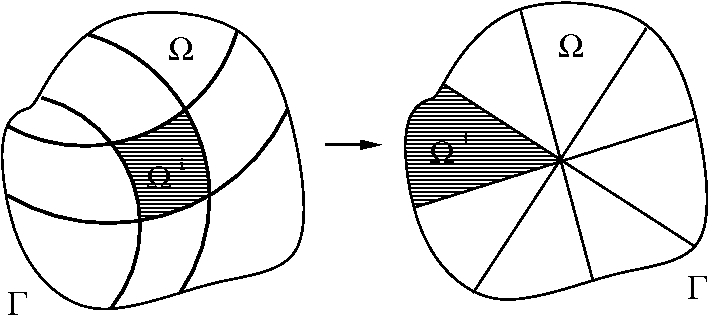}
\caption{Figure presents the domain $\Omega$ and its boundary $\Gamma$. The whole domain $\Omega$ could be divided into subdomains $\Omega^i$ with corresponding line segments $\Gamma^i$ being part of the boundary. The idea of division into subdomains (elements) constitutes the main concept of the finite element method.}
\label{fig:domain}
\end{figure}

The finite elements method (FEM) is based on the idea of division the whole domain $\Omega$ into a number of finite sized elements or subdomains $\Omega^i$ in order to approximate a continuum problem by a behavior of an equivalent assembly of discrete finite elements \cite{zienkiewicz}. 
In the presence of a set of elements $\Omega^i$ the total integral over the domain $\Omega$ is represented by the sum of integrals over individual subdomains $\Omega^i$ 
\begin{equation}
\displaystyle\int_\Omega \mathcal{L}\left(u, \frac{\partial u}{\partial x}, \dots\right)d\Omega = \sum_i\int_{\Omega^i}\mathcal{L}\left(u, \frac{\partial u}{\partial x}, \dots\right) d\Omega ~~{\rm and}~~\int_\Gamma \mathcal{L}\left(u, \frac{\partial u}{\partial x}, \dots\right) d\Gamma = \sum_i\int_{\Gamma^i}\mathcal{L}\left(u, \frac{\partial u}{\partial x}, \dots\right)d\Gamma,
\label{integral}
\end{equation}
where $\mathcal{L}\left(u, \frac{\partial u}{\partial x}, \dots\right)$ denotes a differential operator.
The continuum problem is posed by appropriate differential equations (e. g. Laplace or Poisson equation) and boundary conditions that are imposed on the unknown solution $\phi$. The general procedure of FEM is aimed at finding the approximate solution $\tilde{\phi}$ given by the expansion:
\begin{equation}
\phi \approx \tilde{\phi} = \sum_{j=1}^n \tilde{\phi}^j N_j \overset{\footnotemark}{=} \tilde{\phi}^j N_j,
\label{solution}
\end{equation}    
\footnotetext{the Einstein summation convention}
where $N_j$ are shape functions (basis functions or interpolation functions) \cite{zienkiewicz, kendall} and all or the most of the parameters $\tilde{\phi}_j$ remain unknown. After dividing the domain $\Omega$, the shape functions are defined locally for elements $\Omega^i$. A typical finite element is triangular in shape and thus has three main nodes. It is easy to demonstrate that triangular subdomains fit better to the boundary $\Gamma$ than others e. g. rectangular ones (see Fig.~\ref{fig:domain}). Among the triangular elements family one can find linear, quadratic and cubic elements \cite{zienkiewicz} (see also Appendix \ref{A}). A choice of an appropriate type of subdomains depends on a desired order of approximation and thus arises directly from the continuum problem. The higher order of element, the better approximation. Each triangular element can be described in terms of its \emph{area} coordinates $L^i_{1}, L^i_2$ and $L^i_3$. There are general rules that govern the transformation from \emph{area} to cartesian coordinates
\begin{equation}
\begin{array}{rcl}
x &\!=\!& L_1x_1 + L_2x_2 + L_3x_3 \\
y &\!=\!& L_1y_1 + L_2y_2 + L_3y_3 \\
1 &\!=\!& L_1 + L_2 + L_3
\end{array}
\Leftrightarrow
\left(
\begin{array}{c}
L_1\\
L_2\\
L_3
\end{array}
\right)
=\left(
\begin{array}{ccc}
x_1 & x_2 & x_3 \\
y_1 & y_2 & y_3 \\
1 & 1 & 1
\end{array}
\right)^{-1}
\left(
\begin{array}{c}
x\\
y\\
1
\end{array}
\right)
\label{lcoordiantes}
\end{equation}
where set of pairs $(x_1, y_1), (x_2,y_2), (x_3, y_3)$ represents cartesian nodal coordinates. In turn the \emph{area} coordinates are related to shape functions in a manner that depends on the element order. In further analysis only the linear triangular elements will be used. For them, the shape functions are simply the area coordinates (see Appendix \ref{A}). Therefore, each pair of shape functions $N^i_k(x,y), N^i_l(x,y) {~\rm for}~ k,l=1,2,3$ could be thought as a natural basis of the $\Omega^i$ triangular element.

\subsection{Integral formulas for elements}

We shall consider the linear expression (\ref{linear}) derived in the Appendix \ref{B}
\begin{equation}
\displaystyle \delta\pmb{\tilde\phi}^{T}\left(\mathbf{\pmb{K}}\pmb{\tilde\phi} + \mathbf{f} \right) = \int_\Omega \delta\phi\left\{ -\epsilon \frac{\partial^2 }{\partial x^2} - \epsilon\frac{\partial^2}{\partial y^2} \right\}\phi(x,y) dx dy
    + \int_\Omega\delta\phi\rho dxdy =0
\label{poisson}
\end{equation}
with the boundary condition $\phi=\gamma$ on $\Gamma$. In such a simply case of integral-differential problems with a differential operator $\displaystyle\mathcal{L} = -\epsilon\frac{\partial^2 }{\partial x^2} - \epsilon\frac{\partial^2}{\partial y^2}$, the variable $\pmb{\tilde\phi}$ in Eq.~(\ref{poisson}) only consists of one scalar function $\phi$ which is the sought solution, while the constant vector $\mathbf{f}$ is represented by the last term in the Eq.~(\ref{poisson}). To find the solution for such a problem means to determine the values of $\phi(x,y)$ in the whole domain $\Omega$. The values of $\phi$ on its boundary $\Gamma$ are already prescribed to $\gamma$. On the other hand, at the very beginning (Eq.~(\ref{solution})) we have postulated that a function $\phi$ could be approximated by an expansion $\tilde{\phi}$ given by means of some basis functions $N_m(x,y),~m = 1,\dots,n$ (for more details see Appendix \ref{A}). Thus another possibility to deal with the Poisson problem is just to start from the functional $\Pi$ (Eq.~(\ref{funkcjonal})) and build a set of Euler equations $\displaystyle\frac{\partial \Pi}{\partial \tilde{\phi}_m} = 0$ where $m = 1,\dots,n$ and $\tilde{\phi}_m$ approximates value of the solution $\phi$ calculated at the $m$-th mesh node. 
\begin{equation}
\displaystyle \Pi = \int_\Omega \left\{\frac{1}{2}\epsilon\left(\sum_l\frac{\partial N_l}{\partial x} \tilde{\phi}^l\right)^2 + \frac{1}{2}\epsilon\left(\sum_l\frac{\partial N_l}{\partial y}\tilde{\phi}^l\right)^2 + \rho\sum_l N_l \tilde{\phi}^l\right\}dxdy + \int_\Gamma \left(\gamma - \frac{1}{2}\sum_lN_l\tilde{\phi}^l\right)\sum_k N_k \tilde{\phi}^k d\Gamma
\end{equation}
and after that we calculate the derivative $\displaystyle\frac{\partial \Pi}{\partial \tilde{\phi}_m}$. Moreover, let's simplify our problem by neglecting the last term in the above-presented equation and imposing $\phi = \gamma$ on the boundary $\Gamma$ instead. In that manner, one obtains the expression
\begin{equation}
\displaystyle\frac{\partial \Pi}{\partial \tilde{\phi}_m} = \int_\Omega \left\{\epsilon\left(\sum_l\frac{\partial N_l}{\partial x} \tilde{\phi}^l\right)\sum_k\frac{\partial N_k}{\partial x}\delta^k_{m} + \epsilon\left(\sum_l\frac{\partial N_l}{\partial y}\tilde{\phi}^l\right)\sum_k\frac{\partial N_k}{\partial y}\delta^k_{m} + \rho\sum_l N_l \delta^l_{m}\right\}dxdy = 0
\end{equation}
or in a simplified form
\begin{equation}
\displaystyle\frac{\partial \Pi}{\partial \tilde{\phi}_m} = \sum_l \left\{\int_\Omega \epsilon\left(\frac{\partial N_l}{\partial x}\frac{\partial N^m}{\partial x} + \frac{\partial N_l}{\partial y}\frac{\partial N^m}{\partial y}\right)dxdy\right\}\tilde{\phi}^l  + \int_\Omega\rho N^m dxdy = 0.
\label{continuity}
\end{equation}
It is worth mentioning that some requirements must be imposed on the shape functions $N$. Namely, if $n$-th order derivatives occur in any term of $\mathcal{L}$ then the shape functions have to be such that their $n-1$ derivatives (pay an attention to the Eq.~(\ref{continuity})) are continuous and finite. Therefore, generally speaking $\mathcal{C}_{n-1}$ continuity of shape functions must be preserve.\\
In turn, after substituting
\begin{align}
\vspace{15pt}
\displaystyle K^m_{l} =& \int_\Omega \epsilon\left(\frac{\partial N_l}{\partial x}\frac{\partial N^m}{\partial x} + \frac{\partial N_l}{\partial y}\frac{\partial N^m}{\partial y}\right)dxdy
\label{matrixK}\\
\displaystyle f^m =& \int_\Omega\rho N^m dxdy
\end{align}
finally one obtains a set of equations
\begin{equation}
\displaystyle\frac{\partial \Pi}{\partial \tilde{\phi}_m} = \sum_l K^m_{l}\tilde{\phi}^l  + f^m = 0 ~{\rm for}~ m,l=1\dots,n
\label{important_eq}
\end{equation}
or in matrix description 
\begin{equation}
\pmb{\frac{\partial \Pi}{\partial \tilde{\phi}}} = \mathbf{K}\pmb{\tilde{\phi}}  + \mathbf{f} = 0.
\end{equation}
It is worth noticing that the matrix $\mathbf{K}$ is a symmetric one because of the symmetry in exchange of subscripts $l$ and $m$ in Eq.~(\ref{matrixK}).\\
Now, we are obliged to employ the division of our domain $\Omega$ into a set of subdomains $\Omega^i$. It gives that
\begin{eqnarray}
K^m_l &\!\!=\!\!& \sum_i K^{im}_l \!=\! \sum_i \int_{\Omega^i} \epsilon(x,y)\left(\frac{\partial N^i_l(x,y)}{\partial x}\frac{\partial N^{im}(x,y)}{\partial x} \!+\! \frac{\partial N^i_l(x,y)}{\partial y}\frac{\partial N^{im}(x,y)}{\partial y}\right)dxdy
\label{Kexpression}
\\
f^m &\!\!=\!\!& \sum_i f^{im} \!=\! \sum_i \int_{\Omega^i}\rho(x,y) N^{im}(x,y) dxdy.
\label{function}
\end{eqnarray} 
Therefore, after the transformation to $I$ subdomains the expression (\ref{important_eq}) becomes
\begin{align}
\displaystyle \left( \sum_i K^{im}_{l} \right) \tilde{\phi}^{l}  &+ \left( \sum_i f^{im} \right)  = 0 ~{\rm for}~ i=1,\dots,I~{\rm and}~m=1,\dots,n
\label{important_eq2}\\
\displaystyle \pmb{\tilde{\phi}} &= -\left(\mathbf{K}\right)^{-1} \mathbf{f} \tag{in matrix notation}
\end{align}
In fact, the summation in Eq.~(\ref{important_eq2}) takes into account only these elements $\Omega^i$ which contribute to $m$-th node, however, because of the consistency in notation all elements are included in the sum with the exception that those $N^i_m$ functions for which node $m$ does not occur in $i$-th element are put equal zero.\\  
From now, the whole story is to calculate integrals
\begin{equation}
\begin{array}{c}
\vspace{10pt}
\displaystyle K^{im}_l = \int_{\Omega^i} \epsilon(x,y)\left(\frac{\partial N^i_l(x,y)}{\partial x}\frac{\partial N^{im}(x,y)}{\partial x} \!+\! \frac{\partial N^i_l(x,y)}{\partial y}\frac{\partial N^{im}(x,y)}{\partial y}\right)dxdy \\
\vspace{10pt}
\displaystyle = \int_0^1 dL_1 \int_0^{1-L_1} dL_2 \epsilon(L_1,L_2,L_3) \left |\det J^i\right |\left( \nabla N_{i_l} \mathcal{T}\mathcal{T}^T\nabla^T N^{i_m} \right), \\
\displaystyle f^{im} = \int_{\Omega^i}\rho(x,y) N^{im}(x,y) dxdy = \int_0^1 dL_1 \int_0^{1-L_1} dL_2 \left |\det J^i\right | \rho(L_1,L_2,L_3) N^{i_m}(L_1,L_2,L_3) 
\label{local_integrals}
\end{array}
\end{equation}
where $N^{i_m} = L^{i_m}$, $L_3 = 1 - L_1 - L_2$ (see Appendix \ref{A}) whereas $\det J^i$ -- the jacobian of $i$-th element, $\mathcal{T}$ matrix together with  $\nabla$ operator in new coordinates are evaluated in Appendix \ref{C}.
An integration over the $i$-th subdomain $\Omega^i$, which is a triangular element with three nodes, enforces the transformation from $n$-dimensional global interpolation to the local interpolations given by means of $\displaystyle N_{i_k}(x,y)$ functions where $i_k=1,2,3$. That is why in Eqs (\ref{local_integrals}) new indices $i_l,i_m$ appear which further are allowed to take three possible values 1,2 and 3 for each element $i$ (the local subspace).\\
As a next step, the Gauss quadrature is employed to compute above-written integrals numerically as it is described in the Appendix \ref{D}.
And finally, after incorporating boundary conditions to Eq.~(\ref{important_eq2}) by inserting appropriate boundary values of $\tilde{\phi}$, the system of equations can be solved.

\section{Mesh generation}

\subsection{Initial mesh}

\begin{figure}[t!]
\includegraphics[scale=0.3]{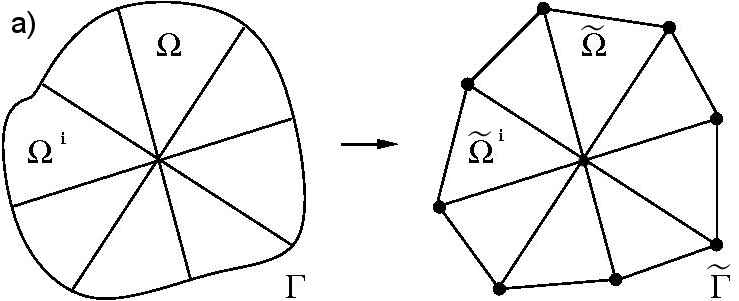}
\includegraphics[scale=0.3]{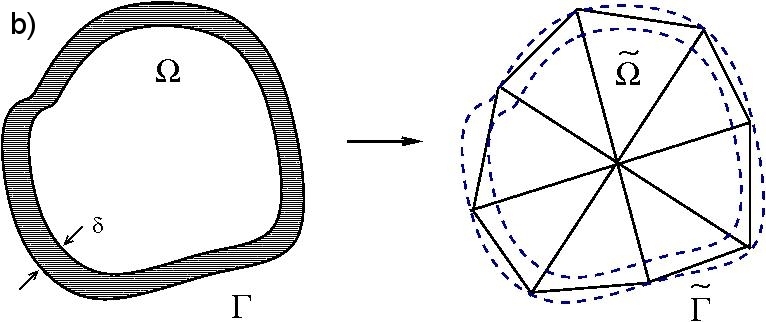}
\caption{Figure presents the domain $\widetilde{\Omega}$ and its boundary $\widetilde{\Gamma}$ after projection to the polygonal domain. It has eight boundary nodes and one central node. Comparing both the initial $\Omega$ and the polygonal $\widetilde\Omega$ domain one can notice that such a simple projection gives rather rough correspondence between them a), however, in some cases it could be a sufficient one i. e. when an integrated function changes very slowly in some $\delta$-thick neighbourhood of the boundary $\Gamma$ b).}
\label{fig:domain2}
\end{figure}

The domain $\Omega$ is a set of points in two-dimensional Euclidean plane $\Re^2$ (see Fig.~\ref{fig:domain}). The initial mesh should define the shape of the domain $\Omega$ or more precisely its boundary $\partial\Omega$. Let's denote the bd($\Omega$) as $\Gamma$. It could form a smooth curve (like a circle) or be a polygon. At the beginning it is necessary to assign the starting set of nodes belonging to $\Gamma$. Taking into account polygon it is obvious that the initial mesh must consist of its vertices, however, in the case of a smooth curve one can choose the initial mesh differently. In the article, the author concentrate on the polygonal domains (see Fig.~\ref{fig:domain2}) that can be formed from a smooth curve after placing some initial nodes on its boundary $\Gamma$ and connecting them by line segments (chords). A role of misplaced boundary nodes will be discussed in Sec. 5. 

Let's start with determining the principal rectangular superdomain as a cartesian product $[x_{min}, x_{max}]\times[y_{min}, y_{max}]$ where $\forall x \in \Omega:~ x \geq x_{min}$ etc. and the following function $mesh_{init}$(vertices, radius) where the variable \emph{vertices} determines the number of its sides and the second one gives the \emph{radius} of its circumscribed circle. For instance, one can make use of the Octave GNU project (free open source) \cite{octave} and creature both the initial nodes $p$ and the initial triangles $t$ arrays in the case of regular polygon of $N$ vertices and lying within a circumscribed circle with a given radius. 

function [p, t] = ${\rm mesh_{init}}$(N, radius)

\quad 1: phi = [0:2*pi/N:2*pi*(N-1)/N]';

\quad 2: p = [radius*cos(phi), radius*sin(phi)]; 

\quad 4: \quad p$^{{\rm center}}$ = sum(p,1)/size(p,1);

\quad 5: p = [p; p$^{{\rm center}}$];

\quad 6: for i = 1:(N-1)

    \quad 7: \quad t(i,:) = [i, i+1, N+1];
    
\quad 8: end

\begin{figure}[htb!]
\includegraphics[scale=0.25]{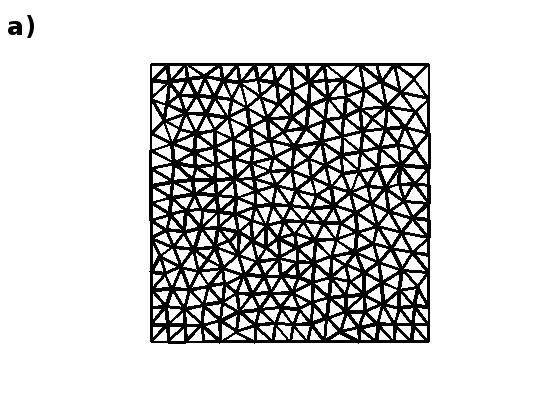}
\includegraphics[scale=0.23]{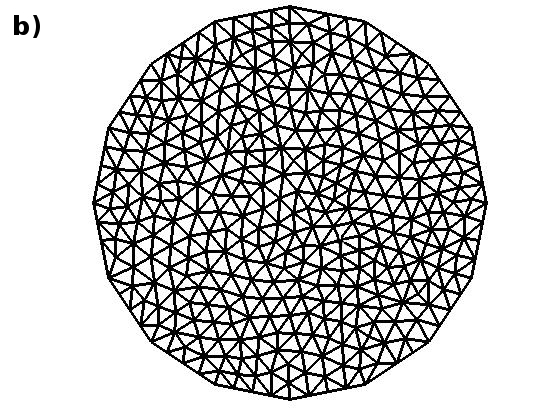}
\includegraphics[scale=0.24]{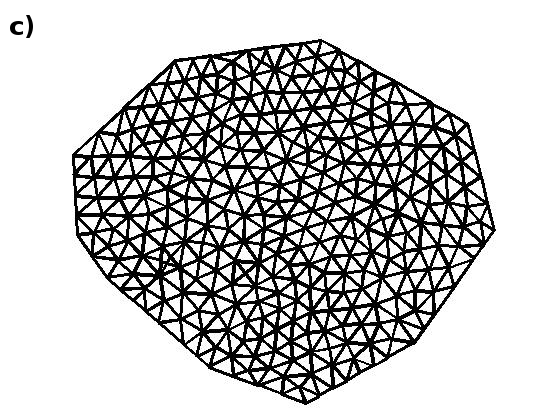}
\includegraphics[scale=0.25]{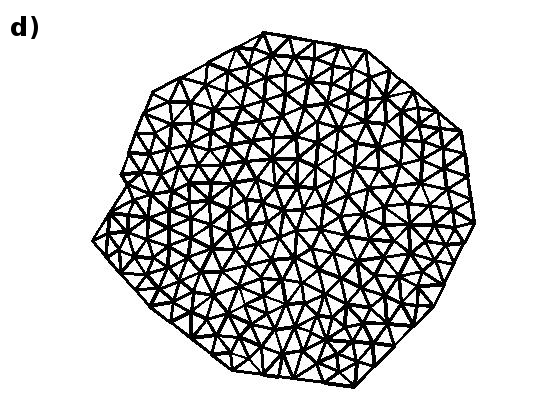}
\includegraphics[scale=0.35]{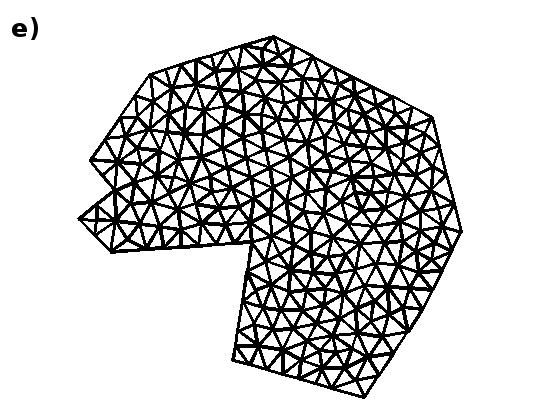}
\includegraphics[scale=0.35]{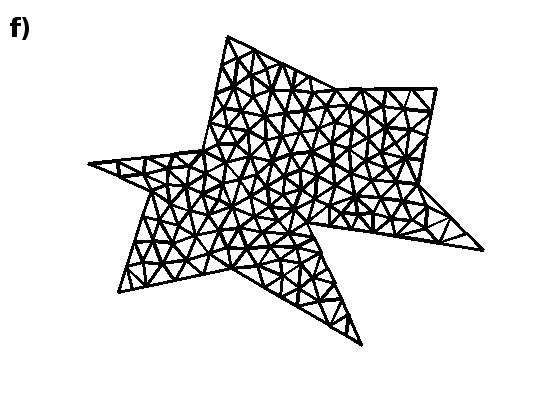}
\caption{Figure shows four domains $\Omega$ having different shapes. In brackets, finally established set of parameters is written: $N_{p}$ -- number of mesh points, $N_{divisions}$ -- number of divisions (according to Sec. 3.2), $\overline{S}_{N}$ -- a normalized average element area are presented; a) regular polygon -- square (258, 8, 1.002); b) regular polygon with 16 nodes (376, 6, 1.026) which approximates circular shape well; c) non-regular, convex figure (315, 8, 1.01); d) non-regular, semi-convex figure (247, 6, 1.071); and two non--regular, non--convex figures e) (245, 7, 0.993) and f) (164, 6, 1.0003) both with weight = [0.25, 0.75].}
\label{fig:domain3}
\end{figure}

\quad 9: t = [t; 1, N, N+1];

end

where $p^{{\rm center}}$ in line 4 denotes the geometrical center of a figure that is an accurate choice for convex cases like regular polygons. However, not only convex type problems are available by presented routine. One can set explicitly $p$ table of initial points and compute $t$ table on its basis (lines 6-9). In such a case the figure center might require to be shifted, for instance, by the formula given in line 13. That center displacement is done in respect to the superdomain center, here set as $[0,~0]$ point. The chosen values of weight vector depend on the problem. In Fig.~\ref{fig:domain3}, meshes for non-convex figures were obtained with weight = $[0.25, 0.75]$.\\

\quad 10: If \emph{non-convex figure}

\quad 11: \quad p1 = p(sum(p.\textasciicircum 2 - repmat(p$^{{\rm center}}$, [size(p,1), 1]).\textasciicircum 2, 2) $>=$ 0, :);

\quad 12: \quad p2 = p(sum(p.\textasciicircum 2 - repmat(p$^{{\rm center}}$, [size(p,1), 1]).\textasciicircum 2, 2) $<$ 0, :);

\quad 13: \quad p$^{{\rm center}}$ = sum([weight(:, 1)*p1; weight(:, 2)*p2], 1)/size([p1; p2], 1);

\quad 14: end

Following further steps of the algorithm presented in next sections, one can obtain meshes for different domains $\Omega$ (see few examples in Fig.~\ref{fig:domain3}).\\
Let's introduce a measure that estimates an element area in respect to the prescribed element area $S$ designed by the element size $h$. The measure $S_{N} = S_{elem}/S$ gives a normalized area for each element. An estimation of the average deviation from assumed value of the element area provides information of mesh quality in the case for their fairly uniform distribution.   

\subsection{Adding new nodes to the mesh}

In this section, let's start with the procedure that allows us to add new mesh nodes to the existing ones.  
The initial configuration of the nodes were already defined. It \textbf{must} define well the shape of the divided area in aspects explained in the description of the Figure~\ref{fig:domain2}. These \emph{initial} nodes are called the \emph{constant nodes} and are kept immobile through the rest of the algorithm steps. Each triangle could be split up into two new triangles by adding a new node to its longest bar. To avoid producing triangles much smaller than defined by the element \emph{size} $h$ only part of them could be broken up i. e. these for which the triangle area is one and half times bigger than $\mathcal{A}$. That condition is set in the algorithm by introducing a new control parameter $C_{split}$. The new node is added in the middle of the triangle longest bar.  

\begin{figure}[htb!]
\includegraphics[scale=0.45]{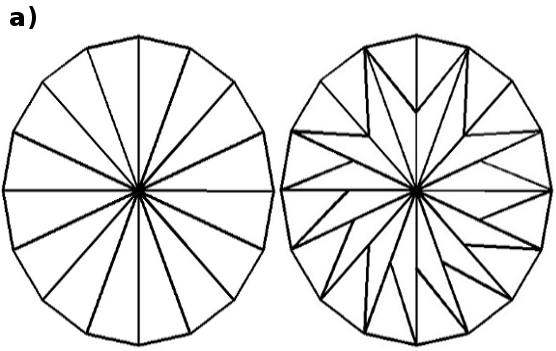}
\hspace{-0.23in}
\includegraphics[scale=0.31]{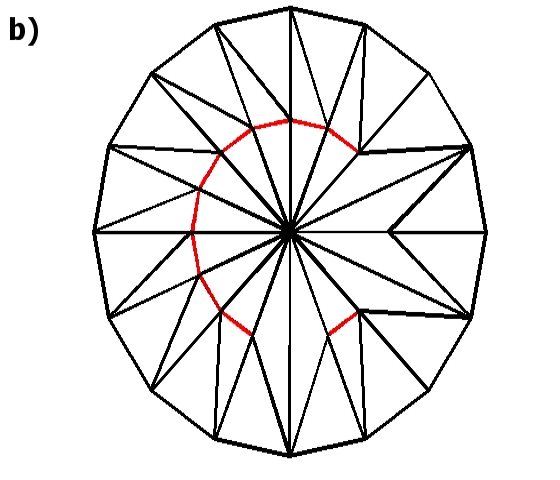}
\includegraphics[scale=0.51]{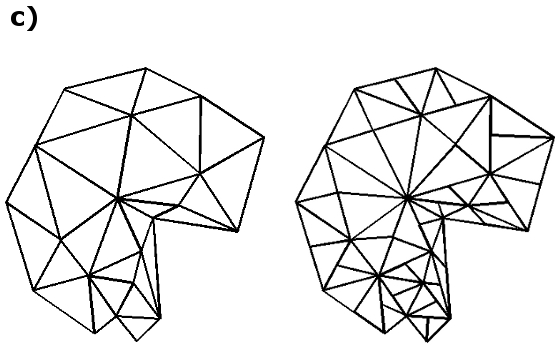}
\hspace{1.2in}
\includegraphics[scale=0.37]{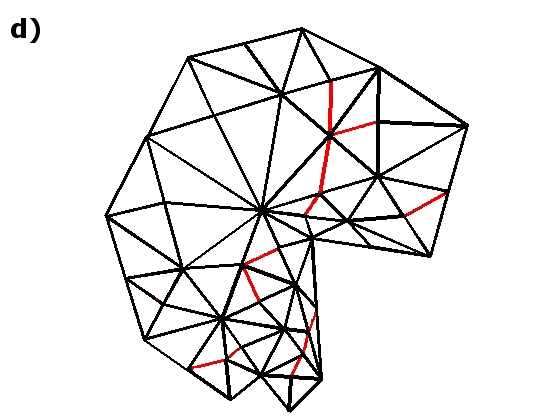}
\caption{Figure presents a division process of non-regular and circular domains together with their boundaries. Pictures a) and c) show meshes with new nodes. Some of them are of the \emph{illegal} type (defined in Sec. 3.2). These nodes constitute starting points for next complementary division that transforms such not well--defined elements into the correct ones, see pictures b) and d).}
\label{fig:domain4}
\end{figure}

For each triangle $\mathcal{T}_k \in \Omega$ for which $C^k_{split} = 1$ 

\quad Find its longest bar $bar^k_{longest}$

\quad Calculate a position of a new node $p^k_{mid}$

\qquad If the node $p^k_{mid}$ is the new one

\qquad Add it to the nodes table

\qquad end

\quad Update triangles table by replacing the old triangle $\mathcal{T}_k$ by two new triangles based on $p^k_{mid}$   
        
end

It is worth mentioning that presented above algorithm is not quite optimal because some of the new nodes could produce triangles with one edge divided by a node resulting from splitting up an adjacent triangle. Such triangles are not desirable and are denoted as $\mathcal{T}_{illegal}$ (see Figures \ref{fig:domain4}a) and c)). Thus the previous procedure needs to be improved. Let's add a few extra steps to it:

For each triangle $\mathcal{T}_k\in\Omega$ perform checking whether it is of $\mathcal{T}_{illegal}$ type

\quad If $\mathcal{T}_k\in\mathcal{T}_{illegal}$

\qquad Split it up into two new properly defined triangles by connecting so-called \emph{illegal} node 

\qquad with the vertex of $\mathcal{T}_k$ lying oppositely to it 

\qquad Remove the old $\mathcal{T}_k$ triangle

\quad end

end

Figures \ref{fig:domain4}b) and d) show meshes having only desired elements.

\subsection{The boundary of the domain}

\begin{figure}[t!]
\includegraphics[scale=0.3]{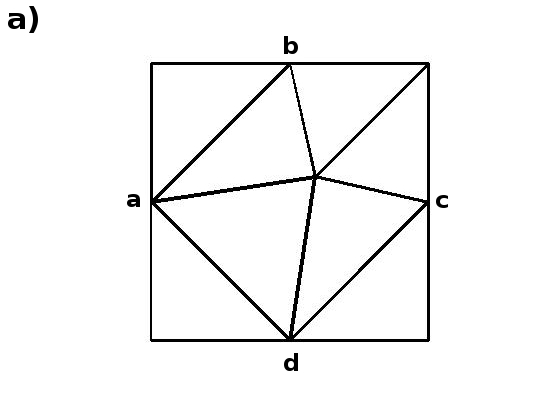}
\includegraphics[scale=0.3]{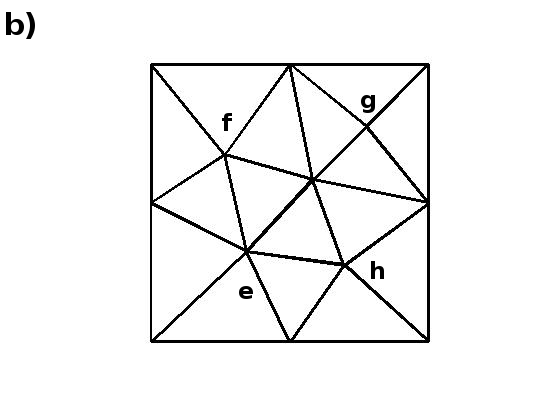}
\caption{
Figure presents the square domain divided into a set of new elements $\widetilde{\Omega}^i$ with corresponding set of line segments $\widetilde{\Gamma}^i$ being its boundary. A way of finding new nodes constitutes the main point of the mesh generation process (see Sec. 3.2) while a selection of nodes is perform according to the algorithm from Sec. 3.3 a) nodes \emph{a,b,c,d} have been classified as boundary nodes whereas b) nodes \emph{e,f,g} have been determined as internal nodes.}
\label{fig:domain5}
\end{figure}

The one of the most important issues to definite is the domain boundary. After determining the boundary $\widetilde{\Gamma}$ by the initial \emph{constant} nodes (lines 1-18 of the presented below algorithm), the next task is to determine which new nodes are lying on boundary line segments $\widetilde{\Gamma}$ (as it is visible in Fig.~\ref{fig:domain5}). These selection is done with a help of the following algorithm:

\quad // For an initial node table $p$ (nodes from 1 to N) find all pairs of neighbouring vertices:

\quad 1: \quad pairs = zeros([ ], 2);

\quad 2: \quad for i = 1:(N-1)

\quad 3: \qquad pairs = [pairs; i i+1];

\quad 4: \quad end

\quad 5: \quad pairs = [pairs; N 1];

\quad // Connect them by a segment line. If $x_1-x_2 \ne 0$ then a function $y=ax+b$ exists and one

\quad // can find pairs $a,b$ for each such a line segment otherwise a vertical line $x=a$ must be found

\quad 6: \quad diff = p(pairs(:,1), :) - p(pairs(:,2), :);

\quad 7: \quad TOL = 1.e-5;

\quad 8: \qquad for i = 1:size(diff, 1)

\quad 9: \qquad\quad     if diff(i, 1) $>$ TOL $||$ diff(i, 1) $<$ -TOL

\quad 10: \qquad\qquad          coeff(i, 1) = diff(i, 2)$./$diff(i, 1);

\quad 11: \qquad\qquad          coeff(i, 2) = p(pairs(i, 1), 2) - p(pairs(i, 1), 1).*a(i);

\quad 12: \qquad\quad      else  

\quad 13: \qquad\qquad          coeff(i, 1) = p(pairs(i, 1), 1);

\quad 14: \qquad\qquad          coeff(i, 2) = [ ] or \emph{a value out of bounds i. e. the principal rectangular superdomain } 

\quad 15: \qquad\quad      end

\quad 16: \qquad\qquad          coeff(i, 3) = min(p(pairs(i,:), 2));

\quad 17: \qquad\qquad          coeff(i, 4) = max(p(pairs(i,:), 2));

\quad 18: \qquad end

Establish the table of coefficients $a,b$ once.

\quad For each new node

\qquad Check whether its coordinates ($x,y$) fulfill any of $y=ax+b$ equations or $x = a$ where $y < y_2$ and $y > y_1$

\qquad\quad If yes classify it as boundary node

\qquad\quad else classify it as internal node

\qquad\quad end

\quad end  

\section{Optimization via the Metropolis method}

Let us define the set of mesh triangles $\Omega = \{ T_j,~~j=1,2,\dots, M \}$ and a set $\mathcal{T}^i$ of triangle mesh elements to which a node $p_i$ belongs. The \emph{closest neighbours} $\mathcal{C}(p_i)$ of the mesh point $p_i$ are defined as a subset of mesh points $p_j \in \mathcal{P}$
\begin{equation}
\forall p_i~~  \exists\mathcal{T}^i \subset \Omega: ~~p_i \in \mathcal{T}^i  ~~ \mathcal{C}(p_i) = \{ p_j:~ p_j \in \mathcal{T}^i ~{\rm for} ~~i\ne j \}.
\end{equation}
Note, that the \emph{closest region} is not the same what \emph{the Voronoi region} \cite{zienkiewicz}. Presented definition is needed to proceed with the Metropolis algorithm \cite{Metropolis} which will be applied in order to adjust triangle's area to the desired size given by the element \emph{size} $h$.
In turn, a proper triangulation is the essence of the finite element method as it is stated in the Sec. 2. Let us divide the whole problem into two different tasks. The first one focuses on finding an optimization for mesh elements being \emph{the internal elements} whereas the second one is developed for so-called \emph{the edge elements}. They are the elements for which one triangle's bar belongs to the boundary $\Gamma$ of the domain $\Omega$. It is assumed that a proper triangulation gives a discrete set of triangles $\mathcal{T}_j$ which approximates the domain $\Omega$ well.\\   
\begin{figure}[htb!]
\includegraphics[scale=0.6]{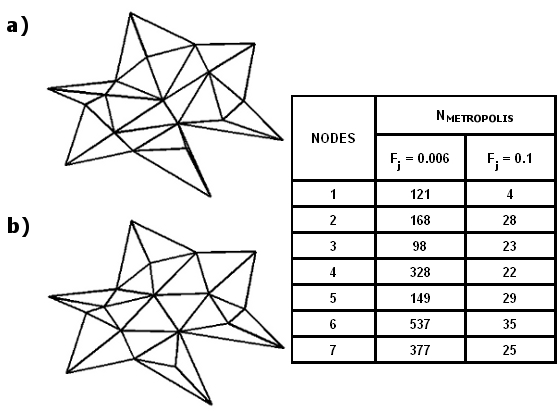}
\caption{Figure shows an application of \emph{the Metropolis algorithm}. Picture a) presents initial positions of new nodes just after generating them whereas picture b) shows their positions after node shifts according to the procedure described in Sec. 4.1 with the following two values of the force strength $F_j$: $0.006$ and $0.1$ applied to each internal node $j=1,\dots, 7$ and temperature set as $T=0.01$. The table presents the total number of Metropolis steps that was required to obtain the final result shown in b).}
\label{fig:metropolis}
\end{figure}

\subsection{The \emph{internal} elements}

Presented method is based on the following algorithm:

\begin{itemize}
\item
Define the element \emph{size} $h$ and consequently the element \emph{area} $\mathcal{A}$.
\item
Initialize the configuration of triangles and then select the \emph{internal} nodes $\mathcal{P}_{int} = \{ p_i: p_i \in \mathcal{P} ~\wedge~ p_i \notin \Gamma \}$ i. e. these nodes does not belong to the domain boundary $\Gamma$.
\item
For each node $p_i$ in $\mathcal{P}_{int}$ find its subdomain $\Omega^{i}$ defined as a set of triangles $\mathcal{T}_i$ to which the node $p_i$ belongs.
\item
Perform the Metropolis approach to every \emph{internal} node $p_i$ within its subdomain $\Omega^{i}$. The Metropolis algorithm is adopted in order to adjust an area of each triangle in the node's subdomain to prescribed value $\mathcal{A}$ by shifting the position of the node $p_i$ (Fig.~\ref{fig:metropolis} demonstrates robustness of the Metropolis approach; compare the node distribution in \textbf{a)} and in \textbf{b)}). That adjustment is govern by the following rules:
\begin{itemize}
\item
Find an area of each triangle $\mathcal{A}_k$ (where $k = 1, 2, \dots, K$) in $\Omega^{i}$ together with the vectors $\vec{r}_{ji} = p_i - p_j$ for each $p_j \in \Omega^i$ connected to node $p_i$
\item
Calculate the length of each triangle edge $\|\vec{r}_{ji}\|$ and its deviation $\delta \|\vec{r}_{ji}\|$ from the designed element size $h$ i. e. $\delta \|\vec{r}_{ji}\| = \|\vec{r}_{ji}\| - h$
\item
Calculate the new position of the node $p^{new}_i$ as
\begin{equation}
\displaystyle p^{new}_i = p_i - \sum_j F_j \delta \|\vec{r}_{ji}\| \frac{\vec{r_{ji}}}{\|\vec{r_{ji}}\|},
\end{equation}
where $F_j$ are weights corresponding to magnitude of $j$-th force applying to node $p_i$. Finally, they were set to the constant value $F$.
\item
Find an area of each triangle $\mathcal{A}^{new}_k$ in $\Omega^{i}$ after shifting $ p_i \to p^{new}_i$
\item
Apply an energetic measure $\mathcal{E}$ to a sub-mesh $\Omega^i$. That quantity could be understand in terms of a square deviation of a mesh element area  from the prescribed element \emph{area} $\mathcal{A}$. Therefore, in the presented paper the $\delta \mathcal{E}$ is defined as a sum of a discrepancy between each triangle \emph{area} $\mathcal{A}_k$ and $\mathcal{A}$ after moving node $p_i$ and prior it, respectively
\begin{equation}
\displaystyle \delta\mathcal{E} = \sum_k\left((\mathcal{A}^{new}_k - \mathcal{A})^2 - (\mathcal{A}_k - \mathcal{A})^2\right)~~k = 1, 2, \dots, K.
\end{equation}
If the obtained value of an \emph{energetic} change is lower than zero the change is accepted. Otherwise, the Metropolis rule is applied i. e. the following condition is checked
\begin{equation}
\displaystyle e^{-\delta E/T} > r
\end{equation}
where $r$ is a uniformly distributed random number on the unit interval $(0, 1)$ and $T$ denotes temperature. 
\item
The above-presented algorithm is repeated unless an assumed tolerance will be achieved.
\end{itemize}
In order to reach a better convergence of the presented method several other improvements could be adopted. For instance, the change in the length of the triangle edge could be an additional measure of mesh approximation goodness. That condition will ensure a lack of elongated mesh elements i. e. elements with very high ratio of its edge lengths (to see such \emph{skinny} elements look at Fig.~\ref{fig:metropolis}a)).    
\end{itemize} 

\subsection{The \emph{boundary} elements}

The Metropolis algorithm applied to boundary nodes slightly differs from the above-described case and could be summarize in the following steps:
\begin{itemize}
\item
Find all the \emph{boundary} or \emph{edge} nodes i. e. nodes for which $p_{k, edge} \in \Gamma$.
\item
Find triangles in the \emph{closest} neighbourhood of the considered $p_{k,edge}$ node. Then calculate an area of each triangle $\mathcal{A}_{l,edge}$.
\item
Calculate the force acting on each boundary node except \emph{the constant nodes} and coming from \textbf{only} two boundary nodes connected to it.
This imposes the following constrain on the motion of the $k$-th node in order to keep it in the boundary $\Gamma$ 
\begin{equation}
\displaystyle F_{k} = -\sum_{j=1}^2 F\delta \|\vec{r}_{jk}\| \frac{\vec{r_{jk}}}{\|\vec{r_{jk}}\|}   
\end{equation}  
where $\delta \|\vec{r}_{jk}\| $ is defined as previously. The force is tangential to the boundary $\Gamma$.
\item
Similarly, find an area of each triangle $\mathcal{A}^{new}_{l,edge}$ after shifting $ p_{k,edge} \to p^{new}_{k,edge}$ according to the force $F_{k}$.
\item
Adopt the Metropolis energetic condition to the boundary case i.e. 
\begin{align}
\displaystyle \delta\mathcal{E} &= \sum_l\left((\mathcal{A}^{new}_{l,edge} - \mathcal{A})^2 - (\mathcal{A}_{l,edge} - \mathcal{A})^2\right)~~l = 1, 2, \dots, L\\
\displaystyle {\rm if} &\qquad e^{-\delta E/T} > r, \tag{accept}\\
&\displaystyle {\rm otherwise}, \tag{reject}
\end{align}
where $T$ denotes temperature and a random number $r \in U(0,1)$ as previously.
\end{itemize} 
The main point of this part is to ensure that the boundary nodes are moved just \textbf{along} the boundary $\Gamma$ (see Appendix \ref{F}).

\section{Results}

Figure \ref{fig:domain6} presents the square domain (with the edge length equal $\sqrt{2}$) initially divided into a set of new elements (the upper picture). Then a mesh generation process can follow two different ways. The first of them, denoted as (1), is done after switching off the Delaunay procedure and leads to the uniform distribution of 512 identical elements with a normalized triangle area $\overline{S}_N = 0.902$ (equal 0.0039). On the other hand, the second way (2) of creating new elements with help of the Delaunay routine gives almost uniform mesh with a normalized average area equal 1.015. Thus employing such an optimization pattern returns a result closer to desired one whereas the (1) way is much faster and in that particular case also does not make use of the Metropolis procedure. That behavior is caused by a symmetry in element and in node distribution, therefore no node shifting is needed. That example should clarify why other method than merely finding the geometrical center of each element is required to enhance a mesh generation routine.   

\begin{figure}[htb!]
\includegraphics[scale=0.5]{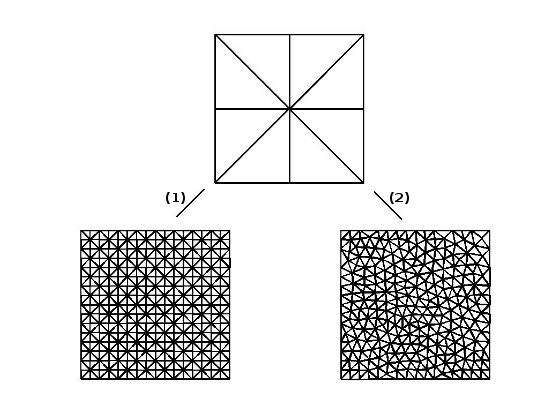}
\caption{
Figure presents two different meshes generated on the basis of a square domain (the upper picture). The way denoted as (1) is done without the Delaunay and the Metropolis optimization procedures giving the uniform distribution of 512 identical elements with $\overline{S}_N = 0.902$. The second way (2) makes use of the Delaunay and the Metropolis routines giving back almost uniform mesh of 459 elements with a normalized average area equal 1.015.}
\label{fig:domain6}
\end{figure}

\begin{center}
{\bf Mesh optimizations}\\
\end{center}

\begin{figure}[htb!]
\includegraphics[scale=0.45, angle=-90]{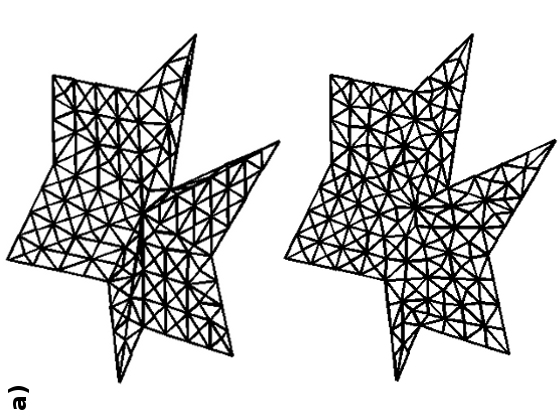}
\includegraphics[scale=0.55, angle=-90]{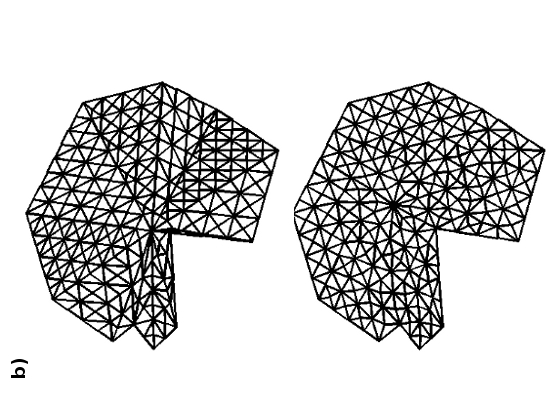}
\caption{
Figure presents two different meshes generated by the algorithm with the Metropolis optimization applied just after new node creation. The obtained mesh characteristics are as they follow: a) before $\overline{S}_N = 1.0589$, $N_p$ = 150, $N_{elem} = 258$, after $\overline{S}_N = 1.0114$, $N_p$ = 160, $N_{elem} = 269$, and b) before $\overline{S}_N = 1.0421$, $N_p$ = 231, $N_{elem} = 416$, after $\overline{S}_N = 1.057$, $N_p$ = 227, $N_{elem} = 406$.}
\label{fig:domain7}
\end{figure}

Figure \ref{fig:domain7} presents results obtained by enriching the proposed mesh generator by the Metropolis approach. Considered meshes were constructed for two non--regular shapes, one of them is also of non--convex type Fig.~\ref{fig:domain7}a). As it is clearly seen in Fig.~\ref{fig:domain7} a mesh quality was in that way enhanced. However, analysis of computed mesh parameters ($\overline{S}_N$ and $N_{elem}$) is not sufficient to explain such mesh improvement. Thus, to quantify meshes with non-uniformly distributed elements (see upper cases of Fig.~\ref{fig:domain7}), the following measure is put forward $\displaystyle S_{var} = mean(\left|S_{elem} - S\right|/S)$. The better fitting to the prescribed element area the smaller value of $S_{var}$. Computed values of $S_{var}$ vary from 0.22 for not optimized meshes to 0.15 after employing the Metropolis rule. Note that for both domains a number of very small elements definitely decreased (see also Fig.~\ref{fig:metropolis}). Moreover, the Metropolis approach offers wide range of feasible mesh manipulations that could be achieved by playing with parameters like the shifting force $F$, the condition of ending optimization (assumed tolerance) and temperature $T$. The shifting force $F$ can differ from node to node or can have the same value for all of them. Furthermore the force strength could change after each node division due to the decline in element areas. The higher value of the accepted force strength $F$, the faster the mesh generation routine. However, putting too high or too law value of $F$ can enormously increase steps of optimization.  

\begin{figure}[htb!]
\includegraphics[scale=0.4]{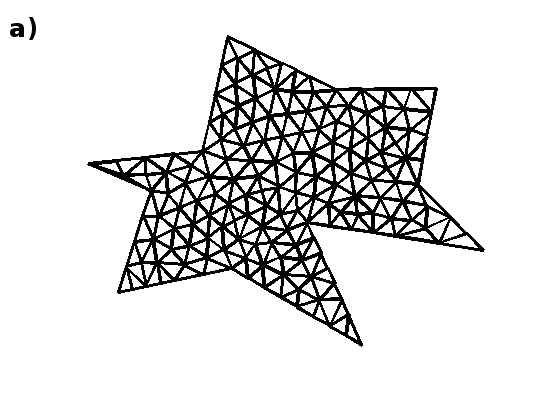}
\includegraphics[scale=0.4]{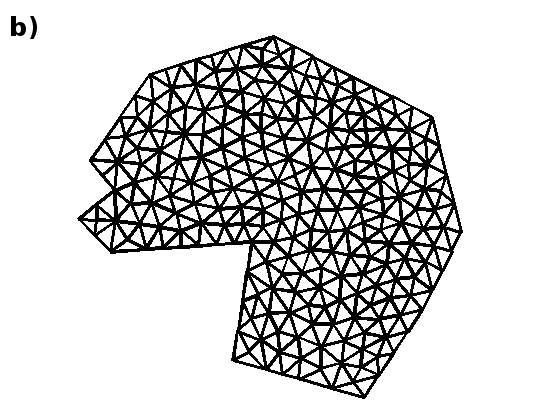}
\caption{
Figure presents two different meshes generated by the algorithm with both the Delaunay and the Metropolis optimizations. The obtained mesh characteristics are as they follow: a) $\overline{S}_N = 0.98$, $N_p$ = 168, $N_{elem} = 278$, and b) $\overline{S}_N = 0.99$, $N_p$ = 244, $N_{elem} = 433$.}
\label{fig:domain8}
\end{figure}

Figure \ref{fig:domain8} shows above--discussed meshes enhanced by additional switching on another kind of optimization i. e. the Delaunay criterion (see Appendix \ref{E}). Both routines were applied in the ordered way i. e.a mesh reconfiguration by the Delaunay method was added before a node shifting via the Metropolis routine. This improves final results in both considered cases. 

\begin{figure}[htb!]
\includegraphics[scale=0.4]{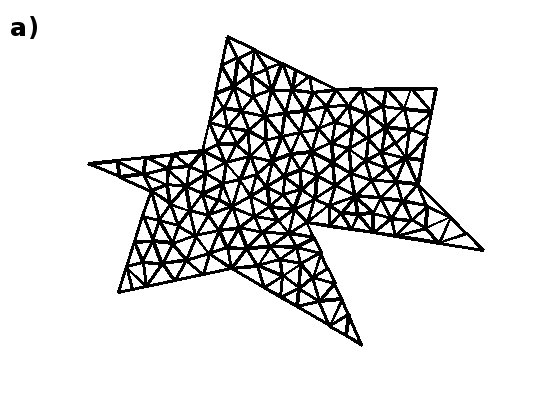}
\includegraphics[scale=0.4]{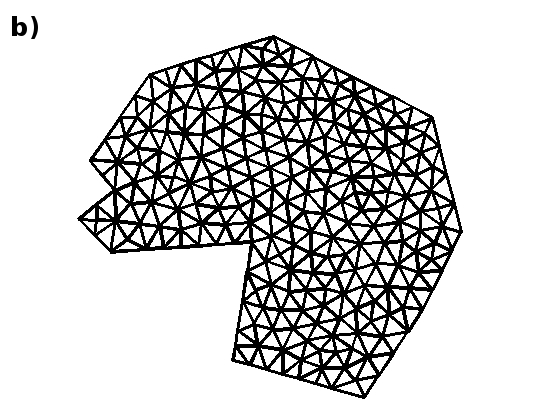}
\caption{
Figure presents two different meshes generated by the main mesh generator (as depicted in previous Figures) but with additional Delaunay routine reconfiguring the old meshes after the Metropolis optimization. New mesh characteristics are as they follow: a) $\overline{S}_N = 1.0003$, $N_p$ = 164, $N_{elem} = 272$, and b) $\overline{S}_N = 0.99$, $N_p$ = 245, $N_{elem} = 432$.}
\label{fig:domain9}
\end{figure}

To examine mesh stability, let's introduce a small perturbation to a mesh configuration obtained by the Metropolis procedure. Applying the Metropolis algorithm leads to the global optimum in element distribution for a given number of nodes, thus resulted mesh should have the stable configuration. To test this, the Delaunay routine was added one more time just after the Metropolis optimization. The highest found changes in mesh quality are depicted in Figure \ref{fig:domain9}. The domains are built with meshes having a little bit different parameters than earlier. In other tested cases no mesh reconfiguration has been detected. The results show that a small exchange in an element configuration in some cases is able to slightly modify the mesh. However, the mesh exchange is hardly seen in the figure \ref{fig:domain9}. Thus, on the other hand, that example demonstrates the stability of the proposed algorithm and proves that the considered mesh generator can be used with confidence. 

\begin{center}
{\bf FEM solutions}\\
\end{center}

\begin{figure}[htb!]
\includegraphics[scale=0.4]{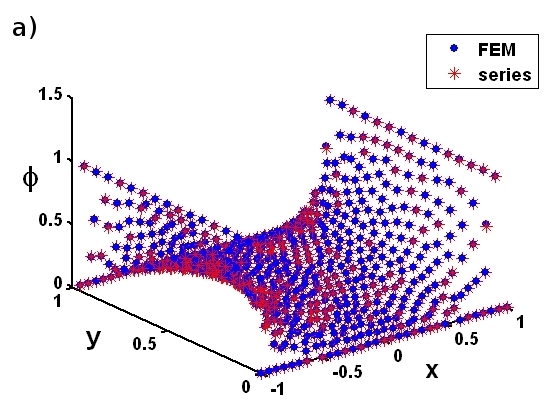}
\includegraphics[scale=0.4]{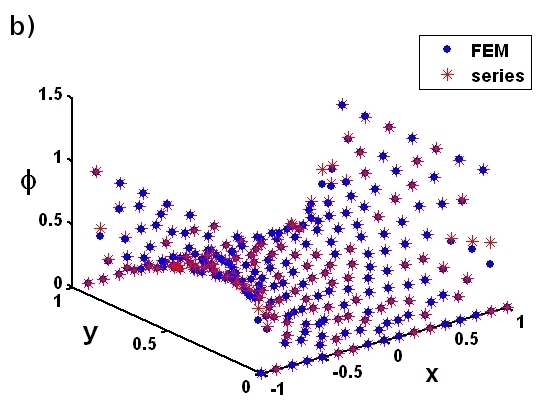}
\caption{
Figure presents two solutions to the Eq.~(\ref{rect_fem}) with the boundary condition $\phi_0$ set as 1. First of them has been obtained by the series (\ref{series}) with $N = 200$ and presents an approximation of the exact solution whereas the second one constitutes a FEM approximation. Computations were performed for two different meshes. The maximum of $\left| \Delta\phi \right| = 0.018$ in the a) case for the element size $h=0.06$ whereas in the case b) max of $\left| \Delta\phi \right| = 0.14$ for the bigger element size $h=0.1$.}
\label{fig:domain10}
\end{figure}

Having above--generated meshes one can perform some computations on them. Thus, let's solve numerically two examples of 2D Poisson problem and then compare them with their exact solutions. The numerical procedure is based on the Finite Element Method already described in Sec. 2. Additionally, in that way mesh accuracy will be tested. The first considered differential problem is embedded within the rectangular domain $[-1 ~1]\times [0 ~1]$ (it constitutes the \emph{mesh}) and has the following form:
\begin{equation}
\displaystyle \frac{\partial^2 \phi}{\partial x^2} + \frac{\partial^2 \phi}{\partial y^2} = 0,
\label{rect_fem}
\end{equation}
with the boundary conditions $\phi=\phi_0$ for $x=-1$ and $x=1$, and $\phi=0$ otherwise.
The solution is given by the series
\begin{equation}
\phi(x,y) = \displaystyle\frac{4\phi_0}{\pi}\sum^N_{n=1,3,5,\dots}\frac{1}{n}\frac{\cosh(n\pi x)}{\cosh(n\pi)}\sin(n\pi y)
\label{series}
\end{equation}
with $N\to\infty$. Figure \ref{fig:domain10} depicts a comparison between an approximation of the analytical solution (Eq.~(\ref{rect_fem})) and a FEM result obtained on the \emph{mesh}. The \emph{mesh} were tested for two $h$ values: 0.1 and 0.06. In the case of a rectangular domain the boundary conditions are fulfilled exactly. Therefore no \emph{boundary} perturbation influences the final result. Analysis of the Fig.~\ref{fig:domain10} confirms a good quality of the mesh allowing to solve accurately the problem under consideration. 

\begin{figure}[htb!]
\includegraphics[scale=0.4]{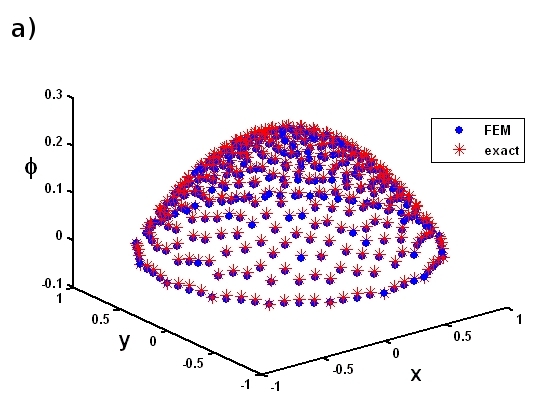}
\includegraphics[scale=0.4]{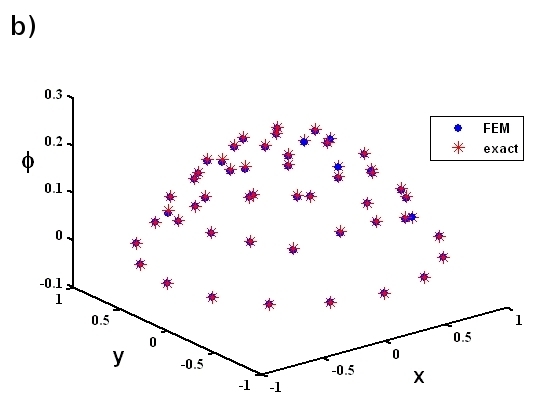}
\caption{
Figure presents comparison between analytical and numerical solutions obtained to Eq.~(\ref{circle_anal}) for two meshes with different element sizes $h$. The maximum of discrepancy between both solutions was also computed and has the following value $\max\left| \Delta\phi \right| = 0.0095$ in the case a) with the element size $h=0.1$ and in the b) case with the element size assumed equal 0.28 $\max\left| \Delta\phi \right| = 0.0067$.}
\label{fig:domain11}
\end{figure}

Let's look on one more differential problem. The following Poisson equation has been solved both numerically and analytically within a circular domain of unit radius  
\begin{equation}
\displaystyle \frac{\partial^2 \phi}{\partial x^2} + \frac{\partial^2 \phi}{\partial y^2} = -1,
\label{circle_anal}
\end{equation}
with the boundary condition $\phi = 0$. The numerical procedure is again based on the FEM approach. The exact solution is given by the expression $\phi = 0.25\left( 1 - x^2 - y^2 \right)$. Figure \ref{fig:domain11} presents the analytical result versus a numerical one. Their comparison shows that discrepancy between them is less than 0.01. Thus, both solutions are in very good quantitative agreement, despite the fact that boundary nodes of the used mesh (Fig.~\ref{fig:domain3}b) do not fulfilled the boundary condition precisely (Fig.~\ref{fig:domain11}a). It is caused by imperfection in this approximation of designed circular shape (see also Fig.~\ref{fig:domain3}). On the other hand, the second mesh (Fig.~\ref{fig:domain11}b) has higher element size ($h=0.28$) than the first one ($h=0.1$) and in that way the boundary condition is fulfilled exactly on the $\widetilde{\Gamma}$. Summarizing, both meshes suffer some weaknesses but they do not influence remarkable solutions.      

All figures presented in the paper were prepared using the Matlab package.

\section{Conclusions}
The proposed mesh generator offers a confident way to creature a pretty uniform mesh built with elements having desired areas. Mesh optimizations are done by means of the Metropolis algorithm and the Delaunay criterion. 
Finding that computed measures $\overline{S}_N$ and $\displaystyle S_{var}$ have pretty good values allows to classify the presented mesh generator as the good one. The meshes were also examined in the context of their stability to some reconfigurations. It was demonstrated that \emph{perturbed} meshes hardly differ from \emph{non--perturbed} ones. Moreover, the mesh was tested by solving the Poisson problem on it making use of the Finite Element Method. The obtained results are in very good quantitative agreement with analytical solutions. This additionally underlines the good quality of the proposed mesh generator as well as efficiency of the FEM approach to deal with differential problems. 

\appendix

\section{The Lagrange polynomials}\label{A}

The Lagrange polynomials $p_k(x)$ are given by the general formula \cite{zienkiewicz, kendall}
\begin{equation}
\displaystyle p_k(x) = \prod^n_{\substack{i=1 \\ i\ne k}}\frac{x-x_i}{x_k-x_i}~{\rm for}~ k=1,\dots,n.
\end{equation}
It is clearly seen from the above-given expression that for $x=x_k~p_k(x_k) =1$ and for $x=x_j$ such that $j\ne k ~p_k(x_j) = 0$. Between nodes values of $p_k(x)$ vary according to the polynomial order i. e. $n-1$ which is the order of interpolation. Making use of these polynomials one can represent an arbitrary function $\phi(x)$ as
\begin{equation}
\phi(x) = \sum_k p_k(x)\phi_k.
\end{equation}
On the other hand, when the interpolated function $\phi$ depends on two spacial coordinates one can define basis polynomials in the form 
\begin{equation}
\displaystyle p_{m}(x,y) \equiv p_{IJ}(x,y) = p_I(x)p_J(y),
\end{equation}
where $I,J$ describe raw and column number for the $m$-th node in a rectangular lattice (rows align along $x$ and columns along $y$ direction, respectively). And consequently, the set $\left\{p_1, \dots, p_m, \dots, p_n \right\}$ is a basis of a $n$-dimensional functional space because each function $p_m$ for $m=1,\dots,n$ equals $1$ at the interpolation node $(x_m,y_m)$ and $0$ at others. It is easy to demonstrate that such functions are orthogonal \cite{kendall}. Instead of spacial coordinates any other coordinates can be considered. In the case of mesh elements the natural coordinates are the \emph{area} coordinates $L$ defined already in the Sec. \emph{The mathematical concept of FEM}. On that basis the shape functions could be constructed as a composition of these three basis polynomials i. e. $N_m(L_1,L_2,L_3) = p^a_I(L_1)p^b_J(L_2)p^c_K(L_3)$ where the values of $a,b$ and $c$ assign the polynomial order in each $L_k$-th coordinate and $I, J$ and $K$ denote the $m$-th node position in a triangular lattice (i. e. in the coordinates $L_1, L_2$ and $L_3$, respectively). In the \cite{zienkiewicz} could be found a comprehensive description of various elements belonging to the triangular family starting from a linear through quadratic to cubic one. For simplicity, in these paper only the linear case is looked on. It explicitly means that shape functions $N_k=L_k(x,y)$, where $k=1,2,3$, change between two nodes linearly (see Eq.~(\ref{lcoordiantes})).

\section{Variational principles}\label{B}

We shall now look on the left hand of the Eq.~(\ref{integral}) i. e. the integral expression $\Pi$
\begin{equation}
\displaystyle \Pi = \int_\Omega \mathcal{L}\left(x, \phi, \frac{\partial \phi}{\partial x}, \frac{\partial \phi}{\partial y}, \dots\right)d\Omega.
\end{equation}
We are aimed at determining the appropriate continuous $\phi$ function for which the first variation $\delta\Pi$ vanish. If   
\begin{equation}
\displaystyle\delta\Pi = \kappa\left(\frac{d}{d\kappa}\Pi[\phi + \kappa\eta]\right)_{\kappa = 0} = 0
\label{var}
\end{equation}
for any $\delta \phi$ then we can say that the expression $\Pi$ is made to be \emph{stationary} \cite{hilbert}.
The function $\phi$ is imbed in a family of functions $\phi +\delta\phi = \phi(x,y) + \kappa\eta(x,y)$ with the parameter $\kappa$. The variational requirement (Eq.~(\ref{var})) gives vanishing of the first variation for any arbitrary $\eta$. In the presented article, the variational problem is limited to the case in which values of desired function $\phi$ at the boundary of the region of integration i. e. at the boundary curve $\Gamma$ are assumed to be prescribed. Generally, the first variation of $\Pi$ has the form 
\begin{equation}
\displaystyle \delta \Pi = \frac{\partial\Pi}{\partial \phi}\delta\phi + \frac{\partial\Pi}{\partial \phi_x}\delta(\phi_x) + \frac{\partial\Pi}{\partial \phi_y}\delta(\phi_y) + \dots
\label{var2}
\end{equation}
and vanishes when
\begin{equation}
\displaystyle  \frac{\partial\Pi}{\partial \phi} = 0, \frac{\partial\Pi}{\partial \phi_x} = 0, \frac{\partial\Pi}{\partial \phi_y} = 0,  \dots
\label{Euler}
\end{equation}
The condition (\ref{Euler}) gives the Euler equations. Moreover, if the functional $\Pi$ is \emph{quadratic} i.e., if all its variables and their derivatives are in the maximum power of 2, then the first variation of $\Pi$ has a standard linear form
\begin{equation}
\displaystyle\delta\Pi \equiv \delta\pmb{\tilde\phi}^{T}\left(\mathbf{K}\pmb{\tilde\phi} + \mathbf{f}\right) = 0,
\label{linear}
\end{equation}
which represents Eq.~(\ref{var2}), though, in matrix notation. A vector $\pmb{\tilde\phi}$ denotes all variational variables i. e. $\phi$ and its derivatives as it is written in Eq.~(\ref{Euler}). $\mathbf{K}$ denotes stiffness matrix (the FEM nomenclature \cite{zienkiewicz}) and $\mathbf{f}$ is a constant vector (does not depend on $\pmb{\tilde\phi}$). We are interested in finding solutions to the Poisson and the Laplace differential equations under some boundary conditions. These classes of differential problems can be represent in such a general linear form (\ref{linear}). Now, we construct a functional $\Pi$ which the first variation gives the Poisson--type equation. Firstly, we define the functional $\Pi$ in the form:
\begin{equation}
\displaystyle \Pi = \int_\Omega \left\{\frac{1}{2}\epsilon\left(\frac{\partial \phi}{\partial x}\right)^2 + \frac{1}{2}\epsilon\left(\frac{\partial \phi}{\partial y}\right)^2 + \rho\phi\right\}dxdy + \int_\Gamma (\gamma - \frac{1}{2}\phi)\phi d\Gamma,
\label{funkcjonal}
\end{equation}
where $\displaystyle d\Gamma = \sqrt{dx^2 + dy^2}$. $\rho,\gamma$ and $\epsilon$ can be functions of spacial variables $x$ and $y$. Secondly, we find the first variation of $\Pi$
\begin{equation}
\displaystyle\delta\Pi = \int_\Omega \left\{ \epsilon \frac{\partial \phi}{\partial x}\delta\phi_x + \epsilon\frac{\partial\phi}{\partial y}\delta\phi_y 
+ \rho \delta\phi \right\}dxdy + \int_\Gamma (\gamma -\phi)\delta\phi d\Gamma
\label{variacja}
\end{equation}
where $\displaystyle\delta\phi_x = \frac{\partial}{\partial x}\delta\phi$. And after integration by parts and taking advantage of the Green's theorem \cite{zienkiewicz} one can simplify the above--written equation to the form
\begin{equation}
\displaystyle\delta\Pi = \int_\Omega \delta\phi\left\{ -\epsilon \frac{\partial^2 \phi}{\partial x^2} - \epsilon\frac{\partial^2\phi}{\partial y^2} 
+ \rho\right\}dxdy + \oint_\Gamma \epsilon \delta\phi\frac{\partial \phi}{\partial n}d\Gamma + \int_\Gamma (\gamma - \phi)\delta\phi d\Gamma= 0,
\label{calkowanie}
\end{equation}
where $\displaystyle\frac{\partial \phi}{\partial n}$ denotes the normal derivative to the boundary $\Gamma$. The expression within the first integral constitutes the Poisson equation
\begin{equation}
\displaystyle -\epsilon \frac{\partial^2 \phi}{\partial x^2} - \epsilon\frac{\partial^2\phi}{\partial y^2} + \rho = 0~~{\rm in}~~\Omega
\end{equation}
whereas the second term in the Eq.~(\ref{calkowanie}) gives the Neumann boundary condition
\begin{equation}
\displaystyle \epsilon \frac{\partial \phi}{\partial n} = 0 {\rm ~~on~~}\Gamma
\end{equation} 
and the third one represents the Dirichlet boundary condition
\begin{equation}
\displaystyle \phi = \gamma {~~\rm on~~}\Gamma.
\end{equation}

\textbf{An important note.} The above--presented calculation demonstrates a way in which one can incorporate the boundary conditions of Neumann or Dirichlet type into a variational expression $\Pi$. However, an appropriately formulated boundary--value problem must include only one kind of b.c. (Neumann or Dirichlet b.c.) defined on the whole boundary $\Gamma$ or is permitted to mix them but in not self--overlapping way. Problems solve in Sec. 5 of the article have only the Dirichlet b.c.. 

\section{Transformation in local $L$-coordinates}\label{C}

Let us compute the determinant of the jacobian transformation between the global $x,y$ and a local $L_1,L_2,L_3$ coordinate frame. One notices immediately that the problem is degenerate. That is why, we introduce a new coordinate $\tilde{z}$ as a linear combination of $L_1,L_2,L_3$ i. e. $\tilde{z} = L_1 + L_2 + L_3$. Note that $\tilde{z}$ is not an independent coordinate and has a constant value equal 1. After taking into account relations (\ref{lcoordiantes}) we find the jacobian matrix in the form
\begin{equation}
\displaystyle J(L_1,L_2,L_3) \equiv \frac{\partial (x,y, \tilde{z})}{\partial (L_1,L_2,L_3)}
= \left(
\begin{array}{ccc}
\vspace{5pt}
\displaystyle \frac{\partial x}{\partial L_1} & \displaystyle\frac{\partial x}{\partial L_2} & \displaystyle\frac{\partial x}{\partial L_3} \\
\vspace{5pt}
\displaystyle \frac{\partial y}{\partial L_1} & \displaystyle\frac{\partial y}{\partial L_2} & \displaystyle\frac{\partial y}{\partial L_3} \\
\displaystyle \frac{\partial \tilde{z}}{\partial L_1} & \displaystyle\frac{\partial \tilde{z}}{\partial L_2} & \displaystyle\frac{\partial \tilde{z}}{\partial L_3} \\
\end{array}
\right) =
\left(
\begin{array}{ccc}
\vspace{5pt}
\displaystyle x_1 & \displaystyle x_2 & \displaystyle x_3 \\
\vspace{5pt}
\displaystyle y_1 & \displaystyle y_2 & \displaystyle y_3 \\
\displaystyle 1 & 1 & 1 \\
\end{array}
\right)
\end{equation}
Furthermore, we have the relation between the jacobian and an element area
\begin{equation}
\det J \equiv 2\Delta, 
\label{jacobian}
\end{equation}
where $\Delta$ denotes the area of a triangle which is based on vertices $(x_1,y_1), (x_2,y_2), (x_3,y_3)$. And finally, we obtain the coordinates transformation rule
\begin{equation}
\displaystyle dxdy = 2 \Delta dL_1dL_2, ~{\rm and}~L_3 = 1 - L_1 - L_2.
\end{equation} 
The relation between the gradient operator $\nabla$ in cartesian and in new coordinates is given by:
\begin{eqnarray}
\displaystyle \left[\frac{\partial}{\partial x}, \frac{\partial}{\partial y}\right] &=& \displaystyle \left[\frac{\partial L_1}{\partial x}\frac{\partial}{\partial L_1} + \frac{\partial L_2}{\partial x}\frac{\partial}{\partial L_2} + \frac{\partial L_3}{\partial x}\frac{\partial}{\partial L_3}, \frac{\partial L_1}{\partial y}\frac{\partial}{\partial L_1} + \frac{\partial L_2}{\partial y}\frac{\partial}{\partial L_2} +   \frac{\partial L_3}{\partial y}\frac{\partial}{\partial L_3} \right] \nonumber \\
\displaystyle &=& \frac{1}{2\Delta}\left[\frac{\partial}{\partial L_1}, \frac{\partial}{\partial L_2}, \frac{\partial}{\partial L_3}\right] 
\left(
\begin{array}{cc}
a_1 & b_1\\ 
a_2 & b_2\\
a_3 & b_3
\end{array}
\right) = \left[\frac{\partial}{\partial L_1}, \frac{\partial}{\partial L_2}, \frac{\partial}{\partial L_3}\right] \mathcal{T}
\end{eqnarray}
where $\displaystyle L_k = (a_kx+b_ky+c_k)/(2\Delta)$ ($k=1,2,3$) and $a_1 = y_2 - y_3, b_1 = x_3 - x_2, c_1 = x_2y_3 - x_3y_2$, the rest of coefficients is obtained by cyclic permutation of indices $1,2$ and $3$.

\section{Numerical integration - the Gauss quadrature}\label{D}

The l.h.s integral $\mathcal{I}$ can be approximated by the $Q$ - point Gauss quadrature \cite{zienkiewicz, qw_1, qw_2, qw_3}
\begin{equation}
\displaystyle \mathcal{I} = \int_0^1 dL_1\int^{1-L_1}_0 dL_2 | \det J| f(L_1,L_2,L_3) \approx \sum^Q_{q=1} f_q(L_1, L_2, L_3)W_q
\label{gauss_qw}
\end{equation}
where $W_q$ denotes the weights for $q$ - points of the numerical integration, and can be found in the Table 5.3 in \cite{zienkiewicz}. 
As it was already said, a set of $N_k(L_1,L_2,L_3)$ shape functions where $k=1,2,3$ can be used to evaluate each $f$ function in the interpolation series which, for instance, in the highest order 10 -- nodal cubic triangular element has the following form
\begin{equation}
\displaystyle f(L_1,L_2,L_3) = \sum^3_{k=1} N_k(L_1,L_2,L_3) f^k + \sum^9_{k=4} N_k(L_1,L_2,L_3)\Delta f^k + N_{10}\Delta\Delta f^{10}
\end{equation}
where $f^k$ are nodal values of $f$ function, $\Delta f^k$ denote departures from linear interpolation for mid-side nodes, and $\Delta\Delta f^{10}$ is departure from both previous orders of approximation for the central node\footnote{it comes from the local Taylor expansion written for finite differences} \cite{zienkiewicz}. For linear triangular elements only the first term is important which gives an approximation
\begin{equation}
f(L_1,L_2,L_3) = \sum_{k=1}^3 L_k f^k.
\end{equation} 
Note, that the r.h.s sum (\ref{gauss_qw}) does not include the jacobian $|\det J |$ that should be incorporated by the weights $W_q$ but it is not (in their values given in Table 5.3 from \cite{zienkiewicz}). Thus let's add the triangle area to the above-recalled formula
\begin{equation}  
|\det J | / 2 \sum^Q_{q=1} f_q(L_1, L_2, L_3)W_q
\end{equation}
and in that way we end up with the final expression for the $Q$ -- point Gauss quadrature.

\section{The Delaunay triangulation algorithm}\label{E}

The author would like to remind briefly the main points of the Delaunay triangulation method \cite{delaunay} together with their numerical implementation using Octave package \cite{octave}. Let $\mathcal{P} = \{ p_i, i = 1,2, \dots N\}$ be set of points in two-dimensional Euclidean plane $\Re^2$. They are called \emph{forming points} of mesh \cite{delaunay}. Let us define the triangle $\mathcal{T}$ as a set of three mesh points 
\begin{equation}
\mathcal{T} = \{ t_j \in \mathcal{P}, j = 1, 2, 3\}.
\label{triangle}
\end{equation}
Using the Delaunay criterion one can generate triangulation where no four points from the set of forming points $\mathcal{P}$ are co-circular:
\begin{equation}
\forall p_i \in \mathcal{P} ~\wedge~ p_i \notin \mathcal{T} \quad \left\| p_i - u \right\| > \rho^2
\label{delaunay_criterion}
\end{equation}
where $u$ is the center of the $\mathcal{T}$ triangle and $\rho$ is its radius.\\ 
The proposed algorithm consists of the following steps:
\begin{itemize}
\item
The triangle's bars are given by the following vectors $\vec{t}_{12}, \vec{t}_{13}, \vec{t}_{23}$ where $\vec{t}_{ij} = \vec{t}_j - \vec{t}_i$, $\vec{t}_i = [t^x_i, t^y_i, 0]$ for each $i \ne j$ and $i,j \in \{1, 2, 3\}$ .
\item
The cross product of each triangle bars defines a plane. The pseudovector $\vec{\mathcal{A}}$ together with its projection on the normal to the plane $\hat{n}$-direction $\mathcal{A}_n$ are found
\begin{align}
\vec{\mathcal{A}} = \vec{t}_{12} \times \vec{t}_{13}\\
\mathcal{A}_n = \hat{n}\circ\vec{\mathcal{A}}
\label{cross}
\end{align}
in order to determine the triangle orientation. If the quantity $\mathcal{A}_n > 0$ the triangle orientation is clockwise, otherwise is counterclockwise. 
\item
The determinant of the square matrix $\mathcal{D(T)}$ is built on the basis of the set of triangle's nodes given by Eq. (\ref{triangle})
\begin{equation}
\mathcal{D(T)} = det\left(
\begin{array}{ccc}
t_1^x & t_1^y & (t_1^x)^2 + (t_1^y)^2 \\
t_2^x & t_2^y & (t_1^x)^2 + (t_3^y)^2 \\
t_3^x & t_3^y & (t_2^x)^2 + (t_3^y)^2 \\
\label{determinant_general}
\end{array}
\right)
\end{equation}
next the following determinant is calculated in order to find out whether a mesh point $p_i$ is outside or inside \emph{the Delaunay zone} (see Fig.~\ref{fig:flipping}) 
\begin{equation}
\mathcal{D(T)}_i = det\left(
\begin{array}{cccc}
t_1^x & t_1^y & (t_1^x)^2 + (t_1^y)^2 & 1 \\
t_2^x & t_2^y & (t_1^x)^2 + (t_3^y)^2 & 1 \\
t_3^x & t_3^y & (t_2^x)^2 + (t_3^y)^2 & 1 \\
p_i^x & p_i^y & (p_i^x)^2 + (p_i^y)^2 & 1 \\
\end{array}
\right)
\end{equation}
for each $p_i \in \mathcal{P} ~\wedge~ p_i \notin \mathcal{T}$. 
\item
If for any point $p_i$ its $\mathcal{D(T)} < 0$ the triangle $\mathcal{T}$ is not the Delaunay triangle (see Fig.~\ref{fig:flipping}a). Then one need to
find other triangles in the closest neighbourhood of the triangle $\mathcal{T}$ corresponding to the number of $p_i$ inside the Delaunay zone and recursively exchange the bars between $\mathcal{T}$ and each of them (see Fig.~\ref{fig:flipping}b).
\begin{figure}[t!]
\includegraphics[scale=0.5]{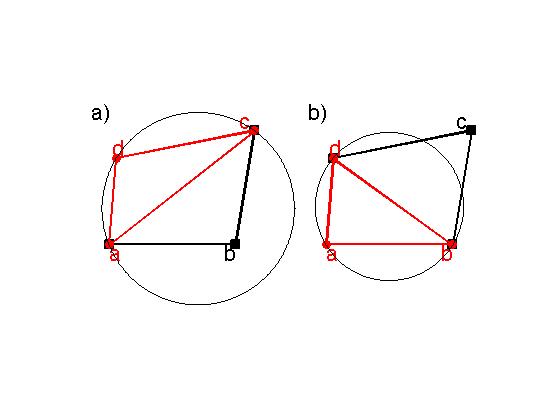}
\caption{Figure shows the main idea of \emph{the Delaunay criterion}. a) Two triangles (with nodes \textbf{abc} and \textbf{acd}, respectively) are not Delaunay triangles, b) after exchange of the edge \textbf{ac} to the edge \textbf{bd} two new triangles \textbf{abd} and \textbf{bcd} replace the old ones. They are both of \emph{the Delaunay type}. Circles represent \emph{the Delaunay zones}.}
\label{fig:flipping}
\end{figure}
\item
Finally, checking whether the new two triangles are the Delaunay triangles takes its place. If so, new ones are accepted unless the change is canceled.
\end{itemize}
Now, let us present the main points of the algorithm that could be easily implemented using the Octave package or Matlab.

Set of mesh points $\mathcal{P} = \{p_i, ~i=1,2,\dots, N \}$ together with their starting configuration forming initial triangular mesh $\Omega = \{\mathcal{T}_j,~ j =1, 2,\dots, M\}$ where $M$ determines \emph{the mesh size}.

while 1
    
    stop = 1
    
    for each mesh triangle $\mathcal{T}$ from $j = 1,2,\dots, M$
    
    \quad ensure the clockwise orientation of the triangle (see above)
    
     \qquad for each mesh point $p_k$ not belonging to the triangle $\mathcal{T}$
     
     \qquad\quad       calculate its $\mathcal{D(T)}_k$
            
     \qquad\qquad       if \textbf{it is lower then 0}
                
     \qquad\qquad\quad           for the triangle $\mathcal{T}$ 
     
     \qquad\qquad\qquad         find its neighbouring $\mathcal{T}^{neighbour}_k$ triangle where the mesh point $p_k$ belongs to
                        
     \qquad\qquad\qquad      exchange the common edge between the triangle $\mathcal{T}$ and its neighbour in order to form 
     
     \qquad\qquad\qquad      two new triangles
     
     \qquad\qquad\qquad      ensure the clockwise orientation of the triangles (see above)
     
     \qquad\qquad\qquad      calculate the newest $\mathcal{D}(\mathcal{T}_{new})_k$ for the new triangle $\mathcal{T}_{new}$
     
     \qquad\qquad\qquad\qquad      if \textbf{ $\mathcal{D}(\mathcal{T}_{new})_k > \mathcal{D(T)}_k$}
     
     \qquad\qquad\qquad\qquad      accept the change and update the set of mesh triangles 
     
     \qquad\qquad\qquad\qquad      put stop = 0
     
     \qquad\qquad\qquad\qquad      else reject the change and return to the previous set of mesh triangles  
     
     \qquad\qquad\qquad\qquad      end
                                    
     \qquad\qquad\quad      end 
                                 
     \qquad\qquad end                     
         
     \qquad  end

     end
     
     if  stop $= 1$ break end

end    

Algorithm ends up with the new triangular mesh $\Omega_{new}$.\\
In order to find an orientation of a triangle $\mathcal{T}$ one can check the sign of $\mathcal{A}_n$ (see Eq.~(\ref{cross})). If it is greater than 0 the triangle orientation is clockwise unless counterclockwise. In the latter case, to ensure the clockwise orientation one can once flip up and down matrix in Eq.~(\ref{determinant_general}) then the triangle orientation turns into the opposite one. Obviously, this flipping results in the change of the sign of the matrix determinant $\mathcal{D(T)} \to -\mathcal{D(T)}$.
    
\section{Boundary nodes -- remarks}\label{F}

For each new boundary node $\tilde{p}$ find its boundary neighbours and save them in the \emph{tboundary} array.
\begin{itemize}
\item       
Find two nodes among all neighbouring nodes in \emph{tboundary} table that are closest to the considered $\tilde{p}$ node. Then save them in \emph{table$_{closest}$}. This works fine for convex figures.     
\item
If the figure's shape is not of convex type then the algorithm must be more sophisticated. It requires to find two nodes that are aligned with the analyzed $\tilde{p}$ node. For example, one can creature an array \emph{table} containing all possible combinations of that node and any other two nodes from the \emph{tboundary} array. Only one combination from the \emph{table} should be the correct one i. e. it must fulfilled conditions describing one of the boundary line segments. Find and save it in \emph{table$_{align}$}.
\end{itemize}

Having \emph{table$_{closest}$} or \emph{table$_{align}$} construct the following \emph{shift$_{dir}$}vector 

\qquad	shift = repmat(-$\tilde{p}$, [2, 1]) + p(table$_{{\rm i}}$, :); i = closest, align

\qquad  shift$_{{\rm dir}}$ = shift./repmat((diag(shift*shift')).\textasciicircum (0.5), [1, 2]);

It defines the node shift direction and must go along with one of boundary line segments.


\begin{thebibliography}{99}
\bibitem{zienkiewicz}
O. C. Zienkiewicz, R. L. Taylor and J. Z. Zhu, \emph{The Finite Element Method: Its Basis and Fundamentals, Sixth edition.}, Elsevier 2005
\bibitem{clough}
R. W. Clough, \emph{The finite element method in plane stress analysis}, In Proc. 2nd ASCE Conf. on Electronic Computation, Pittsburgh, Pa., Sept. 1960; R. W. Clough, \emph{Early history of the finite element method from the view point of a pioneer}, Int. J. Numer. Meth. Eng., \textbf{60}, pp. 283-287, 2004
\bibitem{zienkiewicz2}
O. C. Zienkiewicz, \emph{Origins, milestones and directions of the finite element method}, Arch. Comput. Methods Eng., \textbf{2} (1), pp. 1 - 48, 1995; O. C. Zienkiewicz, \emph{Origins, milestones and directions of the finite element method -- a personal view, Part II: techniques of scientific computing.} In: P.G. Ciarlet and J.L. Lions, editors, \emph{Handbook of Numerical Analysis}, \textbf{4}, pp. 3 - 65, North-Holland, Amsterdam (1996); O. C. Zienkiewicz, \emph{The birth of the finite element method and of computational mechanics}, Int. J. Numer. Meth. Eng., \textbf{60}, pp. 3 - 10, 2004
\bibitem{laplace}
O. C. Zienkiewicz and Y. K. Cheung, \emph{Finite elements in the solution of field problems}, The Engineer, pp. 507-510 1965; O. C. Zienkiewicz, P. Mayer and Y. K. Cheung, \emph{Solution of anisotropic seepage problems by finite elements}, J. Eng. Mech., ASCE, \textbf{92}, pp. 111-120, 1966; O. C. Zienkiewicz, P. L. Arlett, and A. K. Bahrani, \emph{Solution of three--dimensional field problems by the finite element method}, The Engineer, 1967; L. R. Herrmann, \emph{Elastic torsion analysis of irregular shapes}, J. Eng. Mech., ASCE, \textbf{91}, pp. 11-19, 1965; A. M. Winslow, \emph{Numerical solution of the quasi-linear Poisson equation in a non-uniform triangle 'mesh'}, J. Comp. Phys., \textbf{1}, pp. 149-172, 1966;
M. M. Reddi, \emph{Finite element solution of the incompressible lubrication problem}, Trans. Am. Soc. Mech. Eng., \textbf{91}:524 1969 
\bibitem{clough2}
R. W. Clough, \emph{The finite element method in structural mechanics}, In O. C. Zienkiewicz and G. S. Holister, editors, Stress Analysis, Chapter 7. John Wiley \& Sons, Chichester, 1965
\bibitem{delaunay_lit}
A. Bowyer, \emph{Computing Dirichlet tessellations}, Comp. J., \textbf{24}(2), pp. 162-166, 1981;
D. F. Watson, \emph{Computing the n-dimensional Delaunay tessellation with application to Voronoi polytopes}, Comput. J., \textbf{24}(2), pp. 167-172 1981; J. C. Cavendish, D. A. Field and W. H. Frey, \emph{An approach to automatic three-dimensional finite element mesh generation}, Int. J. Numer. Meth. Eng., \textbf{21}, pp. 329-347 1985;N. P. Weatherill, \emph{A method for generating irregular computation grids in multiply connected planar domains}, Int. J. Numer. Meth. Eng., \textbf{8}, pp. 181–197 1988; 
W. J. Schroeder, M. S. Shephard, \emph{Geometry-based fully automatic mesh generation and the Delaunay triangulation}, Int. J. Numer. Meth. Eng., \textbf{26}, pp. 2503-2515 2005; T. J. Baker, \emph{Automatic mesh generation for complex three-dimensional regions using a constrained Delaunay triangulation}, Eng. Comp., \textbf{5}, pp. 161–175 1989; P. L. Georgea, F. Hechta and E. Saltela, \emph{Automatic mesh generator with specific boundary}, Comp. Meth. Appl. Mech. Eng., \textbf{92}, pp. 269-288 1991
\bibitem{afm}
S. H. Lo, \emph{A new mesh generation scheme for arbitrary planar domains}, Int. J. Numer. Meth. Eng., \textbf{21}, pp. 1403-1426 1985; J. Peraire, J. Peiro, L. Formaggia, K. Morgan, O. C. Zienkiewicz, \emph{Finite element Euler computations in three dimensions}, \textbf{26}, pp. 2135-2159 2005; R. L\"{o}hner, P. Parikh, \emph{Three-dimensional grid generation by the advancing front method}, Int. J. Num. Meth. Fluids \textbf{8}, pp. 1135-1149 1988
\bibitem{tm}
M. A. Yerry, M. S. Shephard, \emph{Automatic three-dimensional mesh generation by the modified-octree technique}, Int. J. Numer. Meth. Eng., \textbf{20}, pp. 1965-1990 1984; P. L. Baehmann, S. L. Wittchen, M. S. Shephard, K. R. Grice, and M. A. Yerry, \emph{Robust, geometrically-based, automatic two-dimensional mesh generation}, Int. J. Numer. Meth. Eng., \textbf{24}, pp. 1043-1078 1987; M. S. Shephard and M. K. Georges, \emph{Automatic three-dimensional mesh generation by the finite octree technique}, Int. J. Numer. Meth. Eng., \textbf{32}, pp. 709-749 1991 
\bibitem{phillips}
O. C. Zienkiewicz and D. V. Phillips, \emph{An automatic mesh generation scheme for plane and curved surfaces by isoparametric coordinates}, Int. J. Numer. Meth. Eng., \textbf{3}, pp. 519-528 1971
\bibitem{Peraire}
J. Peraire, M. Vahdati, K. Morgan, and O. C. Zienkiewicz, \emph{Adaptative remeshing for compressible flow computations}, J. Comp. Phys. \textbf{72}, pp. 449-466, 1987
\bibitem{kendall} A. Kendall, H. Weimin, \emph{Theoretical Numerical analysis, A Functional Analysis Framework, Third Edition.}, Springer 2009 
\bibitem{hilbert}
R. Courant, D. Hilbert, \emph{Methods of Mathematical Physics, Volume 1}, Interscience Publisher, New York, 1953
\bibitem{delaunay}
B. Delaunay, \emph{Sur la sph\`{e}re vide}, Izvestia Akademii Nauk SSSR, Otdelenie Matematicheskikh i Estestvennykh Nauk, \textbf{7}, pp. 793-800, 1934
\bibitem{Metropolis}
N. Metropolis, S. Ulam, \emph{The Monte Carlo Method}, J. Amer. Stat. Assoc., \textbf{44}, No. 247., pp. 335-341, 1949
\bibitem{octave}
The source code of Octave is freely distributed GNU project, for more info please go to the following web page	\emph{http://www.gnu.org/software/octave/}.
\bibitem{qw_1}
R. Radau, \emph{\'{E}tude sur les formules d'approximation qui servent \`{a} calculer la valeur d'une int\'{e}grate d\'{e}finie}, Journ. de Math. \textbf{6}(3), pp. 283-336, 1880
\bibitem{qw_2}
P. C. Hammer and O. J. Marlowe and A. H. Stroud, \emph{Numerical Integration Over Simplexes and Cones}, Math. Tables Aids Comp., \textbf{10}, pp. 130-137, 1956
\bibitem{qw_3}
F. R. Cowper, \emph{Gaussian quadrature formulas for triangles}, Int. J. Numer. Meth. Eng., \textbf{7}, pp. 405-408, 1973
\end{thebibliography}
\end{document}